\documentclass[letterpaper, 10 pt, conference]{ieeeconf}

\usepackage{cite}
\usepackage{amsmath,amssymb,amsfonts}
\usepackage{algorithmic}
\usepackage{graphicx}
\usepackage{textcomp}
\usepackage{hyperref}
\usepackage{algorithm}
\usepackage{algorithmic}

\usepackage{graphicx}
\usepackage{epsfig} 
\usepackage{cite}
\usepackage[dvipsnames]{xcolor}

\IEEEoverridecommandlockouts                              

\newtheorem{theorem}{Theorem}[section]

\newtheorem{proposition}[theorem]{Proposition}

\newtheorem{definition}[theorem]{Definition}

\newcommand{\norm}[1]{\left\lVert#1\right\rVert}    

\usepackage{eso-pic}
\AddToShipoutPictureBG*{%
\AtPageUpperLeft{%
\setlength\unitlength{1in}%
\hspace*{\dimexpr0.5\paperwidth\relax}
\makebox(0,-1)[c]{
\begin{tabular}{c c}
Accepted for American Control Conference 2025.
Uploaded to ArXiv \today. \\
\end{tabular}}}}

\AddToShipoutPictureBG*{%
\AtPageUpperLeft{%
\setlength\unitlength{1in}%
\hspace*{\dimexpr0.5\paperwidth\relax}
\makebox(0,-21)[c]{
\footnotesize
\begin{tabular}{c c}
© 2025 IEEE. Personal use of this material is permitted. Permission from
IEEE must be obtained for all other uses, in any \\
current or future media, including reprinting/republishing
this material for advertising or promotional purposes, creating new \\
collective works, for resale or redistribution to servers or lists,
or reuse of any copyrighted component of this work in other works.
\end{tabular}}}}

\begin{document}

\def\BibTeX{{\rm B\kern-.05em{\sc i\kern-.025em b}\kern-.08em
    T\kern-.1667em\lower.7ex\hbox{E}\kern-.125emX}}
\markboth{\journalname, VOL. XX, NO. XX, XXXX 2017}
{Author \MakeLowercase{\textit{et al.}}: Preparation of Papers for IEEE Control Systems Letters (August 2022)}

\title{Online Learning of Interaction Dynamics with Dual Model Predictive Control for Multi-Agent Systems Using Gaussian Processes}

\author{T.M.J.T. Baltussen, E. Lefeber, R. Tóth, W.P.M.H. Heemels, A. Katriniok
\thanks{All authors are with the Eindhoven University of Technology, The Netherlands. Roland Tóth is also with the Systems and Control Laboratory, HUN-REN Institute for Computer Science and Control, Hungary.
{\tt\{t.m.j.t.baltussen, a.a.j.lefeber, r.toth, m.heemels, a.katriniok\}@tue.nl}}}%

\maketitle
\thispagestyle{empty}

\begin{abstract}
The control of a single agent in complex and uncertain multi-agent environments requires careful consideration of the interactions between the agents. 
In this context, this paper proposes a dual model predictive control (MPC) method using Gaussian process (GP) models for multi-agent systems.
While Gaussian process MPC (GP-MPC) has been shown to be effective in predicting the dynamics of other agents, current methods do not consider the influence of the control input on the covariance of the predictions, and hence lack the dual control effect.
Therefore, we propose a dual MPC that directly optimizes the actions of the ego agent, and the belief of the other agents by jointly optimizing their state trajectories as well as the associated covariance while considering their interactions through a GP.
We demonstrate our GP-MPC method in a simulation study on autonomous driving, showing improved prediction quality compared to a baseline stochastic MPC.
The results show that GP-MPC can learn the interactions between the agents online, demonstrating the potential of GPs for dual MPC in uncertain and unseen scenarios.


\end{abstract}

\section{Introduction}
\label{sec:introduction}

Many real-world multi-agent systems are complex and unsuitable for centralized control approaches and require distributed or decentralized controllers \cite{mayne_model_2014}.
Examples include autonomous vehicles (AVs), robotics, aerial vehicles and energy systems, 
which often have coupled dynamics, constraints or policies, making these systems \textit{interactive}.
The control of a single agent requires the awareness and consideration of these interactions with other agents in its environment to achieve safety and performance without relying on communication between the agents.
In the case of AVs, the control of such an agent includes maneuvering and avoiding collisions with other agents.
As AVs often lack explicit communication, we must deduce and consider these uncertain interactions when planning the AV's motion \cite{schwarting_planning_2018}.

Model predictive control (MPC) is an interesting paradigm for this local control as it can directly optimize the agent's trajectories while ensuring safety constraints are satisfied.
Particularly, stochastic MPC (SMPC) can address uncertainties
by incorporating a probabilistic description of model uncertainty into 
a stochastic optimal control problem (SOCP).
In SMPC, the control inputs not only affect the system state, but also the probability distribution associated with this state, which is referred to as the dual control effect \cite{bar-shalom_dual_1974, mesbah_stochastic_2018}.
This dual control effect is two-fold. Firstly, control policies with the dual control effect causally anticipate future observations, and, hence, will be less conservative in regions of low uncertainty. Secondly, dual control policies can feature an inherent probing effect that facilitates \textit{active learning} \cite{bar-shalom_dual_1974}.
Although solving the SOCP via the Bellman equation naturally yields these aspects of the dual control effect, this is often intractable and requires some approximation which may eliminate the dual control effect \cite{mesbah_stochastic_2018}.
In this work, we propose a learning-based MPC method aimed towards advancing dual MPC in complex and dynamic environments.

Learning-based MPC uses data-driven modeling and machine learning methods to improve the system model and the parameterization of the SOCP.
In particular, Gaussian processes (GPs) are an often-used modeling approach for non-parametric learning-based MPC,
as they can compensate for complex model mismatch and enable direct assessment of the approximate model uncertainty \cite{hewing_learning-based_2020}.
As GPs enable universal and differentiable regression \cite{rasmussen_gaussian_2006}, joint conditioning of states \cite{trautman_unfreezing_2010}, and tractable uncertainty propagation \cite{hewing_cautious_2020}, they have strong potential for interaction-aware dual MPC.

\subsection{Related Work}
Dual MPC is shown to be effective in simultaneous identification and control of systems with parametric and structural uncertainties \cite{arcari_dual_2020}.
In addition, sampling-based methods have recently demonstrated the potential of interaction-aware control \cite{wang_interaction_2023}, and active learning in multi-agent environments \cite{knaup_active_2024}.
However, these methods currently rely on parametric models with fixed model structures. By extending these methods to non-parametric models, a larger class of problems can potentially be addressed.
Gaussian process MPC has been used to learn and predict the motion of other agents in autonomous racing \cite{brudigam_gaussian_2021, zhu_gaussian_2023}, and autonomous driving \cite{bethge_model_2023}.
However, previous works fix the predictions of the other agents \cite{brudigam_gaussian_2021} or the covariance of these predictions \cite{zhu_gaussian_2023, brudigam_gaussian_2021, bethge_model_2023} prior to solving the MPC problem, making the controller agnostic to the interactions. Consequently, these methods lack the dual control effect.
While \cite{zhu_gaussian_2023, bethge_model_2023} rely on offline training data, using online training data can help to reduce uncertainty and to improve the generalizability of learning-based control in unseen scenarios, as can be seen in \cite{brudigam_gaussian_2021}.

\subsection{Contributions}
In this paper, we present a novel learning-based dual MPC method to
learn the unknown dynamics and/or policy of other agents
and bridge the previously identified research gap 
through the following contributions:
\text{(i) We} extend GP-MPC for the local control of a single agent in a multi-agent environment by including the dual control effect through a custom GP implementation.
Although we do not incentivize active learning, our MPC explicitly considers the joint covariance of the agents in its optimization.
The MPC leverages the joint probability of the GP to adapt the state covariance of the other agents through the dual control effect, using Bayesian inference to predict the interactions between the ego agent and the other agents.
(ii) Our GP-MPC uses online measurements to update the GP. Through online learning, the GP-MPC can adapt its predictions and handle uncertain and unseen behavior of other agents.
\text{(iii) We} apply the proposed method in a simulation study on lane merging with an autonomous vehicle.
In this study, the GP-MPC learns the closed-loop policy of another agent and adapts the constraints based on its confidence.
To the authors' best knowledge, this is the first work that uses GPs to learn and predict the interactions between agents online via dual MPC.

\section{Problem Formulation}
\label{sec:Problem}
We consider a multi-agent system of {\small$n_a \in \mathbb{N}_{\geq 2}$} agents and the task of devising a control policy for the ego \text{Agent $0$}, which is described by a discrete-time dynamical system

\vspace{-3pt}
\small
\begin{equation}
\label{eq:Ego_Dyn}
    x^0_{k+1} = f^0(x^0_k, u^0_k)
    \vspace{-2pt}
\end{equation}
\normalsize
that is subject to state and input constraints, {\small $x^0_k \in \mathbb{X}^0 \subseteq \mathbb{R}^{n_x}$, $u^0_k \in \mathbb{U}^0 \subseteq \mathbb{R}^{n_u}$} at discrete time step {\small$k \in \mathbb{N}$}.
We assume to have perfect knowledge of the function {\small$f^0$}.
The dynamics of \text{Agent {\small$j \in \mathcal{J}$}}, where {\small$\mathcal{J} = \{1,2,\dots, n_a-1\}$}, are composed of a known, nominal function {\small$f^j$}, and, an unknown and uncertain, residual function {\small$g^j$} that influences the state via the full column rank matrix {\small$B^j \in \mathbb{R}^{n_x \times n_g}$}:

\vspace{-3pt}
\small
\begin{equation}
    \label{eq:Target_Dyn}
    x^j_{k+1} = f^j\bigl(x^j_{k}\bigr) + B^j g^j \bigl( \mathbf{x}_{k}, \mathbf{u}_{k} \bigr) + w^j_k,
    \vspace{-2pt}
\end{equation}
\normalsize
where {\small${\mathbf{x}_{k} = \mathrm{vec}\bigl(x_{k}^0, x_{k}^1, \dots, x_{k}^{n_a-1} \bigr)}$} is the joint state vector, and {\small$\mathbf{u}_{k} =  \mathrm{vec}\bigl(u_{k}^0, u_{k}^1, \dots, u_{k}^{n_a-1} \bigr)$} is the joint input vector, and {\small$w^j_k \in \mathbb{R}^{n_x}$} is a Gaussian disturbance process.
Note that the residual dynamics {\small$g^j$} 
depends on {\small$\mathbf{x}_k$} and {\small$\mathbf{u}_k$}, making the agents \textit{interactive}.
We assume that the control input of Agent {\small$j \in \mathcal{J}$} is composed of a closed-loop state feedback policy:

\vspace{-5pt}
\small
\begin{equation}
    \label{eq:Target_Pol}
    u^j_k= \kappa^j (\mathbf{x}_k).
    \vspace{-2pt}
\end{equation}
\normalsize
As Agent 0 shares its environment with Agent {\small$j \in \mathcal{J}$}, we have to consider their coupled safety constraints that define a safe set
{\small $ {\mathbb{X}^j_k\bigl(x^j_k\bigr) = \{x^0_k \in \mathbb{R}^{n_x} \, \mid \, h^j\bigl(x^0_k, x^j_k  \bigr) \leq 0  \}}$}.
We assume that $f$, $g$, $h$ and $\kappa$ are continuously differentiable functions.
We consider the problem of finding a controller that safely controls \text{Agent $0$} in the sense that {\small$ x^0_k \in \bigcap_{j \in \mathcal{J}} \mathbb{X}^j_k\bigl(x^j_k \bigr) \cap \mathbb{X}^0$}.

\section{Dual Stochastic Model Predictive Control}
\label{sec:Method}
\subsection{Modeling Interactions Between Agents}
In this section, we propose a learning-based dual MPC for the control of Agent $0$.
Since we do not have full knowledge of the dynamics of Agent {\small$j \in \mathcal{J}$}, nor of its policy, i.e., of its intentions,
it is essential that we consider the uncertainty of the residual dynamics {\small$g^j$} and their policy {\small{$\kappa^j$}} in the formulation of our MPC formulation in order to realize the safety of the controller as these interactions will affect the safe set {\small$\mathbb{X}^j_k(x^j_k)$}.
As we only have partial information of the other agents, we predict their future states through Bayesian inference using GP regression.
Firstly, we collect measurements of the states of all agents and the ego's input {\small $\mathbf{z}_k = [\mathbf{x}_k^\top, u_k^{0\top}]^\top \in \mathbb{R}^{n_z}$}, and we collect samples of the residual dynamics {\small$y^j_k \in \mathbb{R}^{n_g}$} as follows:

\vspace{-8pt}
\small
\begin{equation}
    y^j_k = B^{j^\dag} \bigl( x_{k+1}^j - f^j(x_{k}) \bigr) = g^j\bigl(\mathbf{x}_k, \kappa^j(\mathbf{x}_k), u^0_k\bigr) + \Tilde{w}^j_k,
\vspace{-1pt}
\end{equation}
\normalsize
with {\small$B^{j^\dag} \hspace{-1.5mm}=\hspace{-0.5mm} \bigl(B^{j\top} B^j\bigr)^{-1}B^{j\top}$} and {\small$\Tilde{w}^j_k \sim \mathcal{N}(0, \Sigma^j_w) $}.
Then, we use the collected training data 
\text{\small$\mathbf{Z} = [\mathbf{z}_1, \mathbf{z}_2, \dots, \mathbf{z}_{n_D}] \in \mathbb{R}^{n_z \times n_D}$}, {\small${\mathbf{y}^j = [y_1^j, y_2^j, \dots, y_{n_D}^j ] \in \mathbb{R}^{n_g \times n_D}}$}
in order to approximate the residual dynamics by a GP {\small$\mathbf{d}^j \hspace{-0.5mm} : \mathbb{R}^{n_z} \rightarrow \mathbb{R}^{n_g}$} such that:

\vspace{-4pt}
\small
\begin{equation}
    g^j\bigl(\mathbf{x}_k, \kappa^j(\mathbf{x}_k), u^0_k\bigr) 
    \approx \mathbf{d}^j\left( \mathbf{z}_k \right).
    \vspace{-2pt}
\end{equation}
\normalsize

In order to predict the future states of the other agents, we condition the posterior of the GP on the training data {\small$\mathcal{D} = \bigl\{\mathbf{Z}, \bigl( \mathbf{y}^j \bigr)_{j=1}^{n_a - 1} \bigr\}$}.
For the sake of simplicity, we confine ourselves to scalar GPs {\small$\left(n_g =  
 1\right)$}.
Typically, multiple variables are predicted with {\small$n_g \leq n_x$} independent GPs. We refer the interested reader to \cite{hewing_cautious_2020}, for details on vectorial GPs.

\subsection{Gaussian Process Regression}
\label{sec:GP-Reg}
Let us focus on predicting the dynamics of a single \text{Agent {\small$j \ \in \mathcal{J}$}} and omit the superscript {\small$j$} for the sake of readability.
Firstly, we impose a prior distribution on the GP {\small$d\left(\mathbf{z}\right) \in \mathbb{R}$} through a user-defined kernel function {\small$k\left(\mathbf{z}, \mathbf{z}'\right)$} \cite{rasmussen_gaussian_2006}. As we use the GP for model augmentation, we employ a zero mean prior.
We assume that the measured training outputs {\small$\mathbf{y}$}, and the prior distribution of the GP {\small$d\left(\mathbf{z}\right)$} are jointly Gaussian:

\vspace{-2pt}
\small
\begin{equation}
    \begin{bmatrix}
        \mathbf{y} \\ d \left(\mathbf{z}\right)
    \end{bmatrix} \sim \mathcal{N} \left(
        \mathbf{0}, \begin{bmatrix}
            K_{\mathbf{Z Z}} 
            + I \sigma_w^2
            & \mathbf{k}_{\mathbf{Z z}} \\ \mathbf{k}_{\mathbf{z Z}} & k_{\mathbf{z z}}s^2
        \end{bmatrix} \right),
\end{equation}
\normalsize
where {\small$k_{\mathbf{z z}'} = k\left(\mathbf{z}, \mathbf{z}'\right)$}, {\small$\mathbf{k}_{\mathbf{Z z}} \in \mathbb{R}^{{n_D}}$} is the concatenation of the kernel function evaluated at the test point {\small$\mathbf{z}$} and the training set {\small$\mathbf{Z}$}, where {\small$[\mathbf{k}_{\mathbf{Z z}}]_i = k\left(\mathbf{z}_i,\mathbf{z}\right)$} and {\small${\mathbf{k}_{\mathbf{Z z}}^\top = \mathbf{k}_{\mathbf{z Z}}}$}, and {\small$K_\mathbf{Z Z} \in \mathbb{R}^{{n_D} \times {n_D}}$} is a Gram matrix, which satisfies {\small${[K_{\mathbf{Z Z'}}]_{ij} = k\left(\mathbf{z}_i,\mathbf{z}'_{j}\right)}$}.
The Gaussian noise on {\small$\mathbf{y} \in \mathbb{R}^{1 \times n_D}$} is typically handled by regularizing the Gram matrix with the noise variance $\sigma_w^2$, while the predictions on $d$ are noise free.
Marginalizing the jointly Gaussian prior over the training outputs {\small$\mathbf{y}$} yields the GP posterior \cite{rasmussen_gaussian_2006}, which is our conditioned belief of {\small$d(\mathbf{z})$} due to the observations {\small$\mathcal{D}$}:

\vspace{-3pt}
\small
\begin{equation}
    \text{Pr}\bigl(d(\mathbf{z}) \mid \mathbf{z}, \mathcal{D}\bigr) = \mathcal{N}\bigl(\mu^{d}\left( \mathbf{z} \right), \Sigma^{d}\left( \mathbf{z} \right) \bigr),
    \vspace{-2pt}
\end{equation}
\normalsize
where the posterior mean and covariance functions read as:

\vspace{-8pt}
\small
\begin{subequations}
\begin{align}
    \mu^{d}(\mathbf{z}) &= \mathbf{k}_{\mathbf{z Z}} \left(K_{\mathbf{ZZ}} + I \sigma^2_w \right)^{-1} \mathbf{y}, \label{eq:Sparse_Mean}\\
    \Sigma^d(\mathbf{z}) &= k_{\mathbf{z z}} - \mathbf{k}_{\mathbf{z Z}} \left(K_{\mathbf{ZZ}} + I \sigma^2_w \right)^{-1}  \mathbf{k}_{\mathbf{Z}\mathbf{z}} \label{eq:Sparse_Cov}.
\end{align}
\vspace{-2pt}
\label{eq:Sparse_GP}%
\end{subequations}
\normalsize
Subsequently, we query the posterior of the GP {\small$d^j(\mathbf{z})$},
at a test point {\small$\mathbf{z} \in \mathbb{R}^{n_z}$},
to predict the residual dynamics {\small$g^j$} of Agent {\small$j$} over a time prediction horizon {\small$i = 0, 1,\dots, N-1$}:

\vspace{-8pt}
\small 
\begin{subequations} 
    \begin{align}
    x^j_{i+1 \mid  k} & = f^j\bigl(x^j_{i \mid  k}\bigr) + B^j d^j\left( \mathbf{z}_{i \mid k} \right) , \\
    d^j\left(\mathbf{z} \right) & \sim \mathcal{N} \bigl(\mu^{d^j}\left(\mathbf{z} \right), \Sigma^{d^j}\left(\mathbf{z} \right) \bigr),
    \end{align}
    \vspace{-1pt}
\end{subequations}
\normalsize
where {\small$x_{i \mid k}$} denotes the {\small$i$}-th step-ahead prediction of {\small$x_{k+i}$} made at time step $k$ given the data {\small$\mathcal{D}_k$}, where {\small$\mu^{d^j}$} and {\small$\Sigma^{d^j}$} denote the predicted mean and covariance function of the posterior distribution of the GP at a given test point {\small$\mathbf{z}$}.

\subsection{Gaussian Process Prediction Model}
\label{sec:GP-Model}%
We exploit the joint probability of the GP {\small$d^j$} to account for the uncertainty in our safety constraint set {\small$\mathbb{X}^j$} by formulating \textit{chance constraints} that should be satisfied with a probability \text{\small $p_x \in \left(0, 1\right]$}.
To this end, we optimize the control sequence of Agent 0 \text{{\small $U_k=\bigl( u^0_{0 \mid k}, \ldots, u^0_{N-1 \mid k}\bigr)$}} in the following SOCP:

\vspace{-8pt}
\small
\begin{subequations}
\label{eq:OCP}
\begin{align}
\min _{U_k} & \quad  J \bigl( x^0_k, U_k \bigr) \\
\text {s.t. } x^0_{i+1 \mid k} & = f^0 \bigl(x^0_{i \mid k}, u^0_{i \mid k}\bigr), \quad i = 0,1, \dots , N-1, \\
x^j_{i+1 \mid k} & = f^j \bigl( x^j_{i \mid k} \bigr) + B^jd^j\bigl( \mathbf{z}_{i \mid k} \bigr) \notag \\
i  & = 0, 1, \dots, N-1, \quad j \in \mathcal{J}, \\
\text{Pr} \Bigl(  x^0_{i \mid k} \in \mathbb{X}^j_{i \mid k} &\bigl(x^j_{i \mid k}\bigr) \Bigr) \geq p_x, \quad i = 1,2, \dots, N, \; \; j \in \mathcal{J} \label{eq:OCP_CC}\\
x^0_{i \mid k} & \in \mathbb{X}^0, \quad i = 1,2, \dots , N, \\
u^0_{i \mid k} & \in \mathbb{U}^0, \quad i = 0,1, \dots , N-1, \\
x^j_{0 \mid k} & = x^j_k, \quad j \in \{0 \cup \mathcal{J}\}.
\end{align}
\end{subequations}
\normalsize
Note that the posterior GP {\small$d^j(\mathbf{z})$} is a function of the control sequence $U_k$.
As SOCP \eqref{eq:OCP} is generally intractable, we subsequently apply tractable approximations for consecutive evaluations of the GP and the chance constraints at the expense of formal safety guarantees.
In order to propagate the covariance of the predicted state {\small$\Sigma^{x^j}$}, we make the following assumptions. 
We assume that consecutive GP evaluations are independent and the predicted state {\small$x^j$} and posterior GP {\small$d^j(\mathbf{z})$} are assumed to be jointly Gaussian at each prediction step:

\vspace{-8pt}
\small
\begin{equation}
\label{eq:Joint_Dist}
\begin{bmatrix}
    x^j_{i} \\
    d^j(\mathbf{z}_i)   
\end{bmatrix}
 \sim \mathcal{N}\left(\begin{bmatrix}
\mu^{x^j}_{i} \\
\mu^{d^j}(\mathbf{z}_i)
\end{bmatrix},
\begin{bmatrix}
\Sigma^{x^j}_{i} & \hspace{-2mm} \Sigma^{x^j d^j}(\mathbf{z}_i) \\
\Sigma^{d^j x^j}(\mathbf{z}_i) & \hspace{-2mm} \Sigma^{d^j}(\mathbf{z}_i)
\end{bmatrix}\right),
\end{equation}
\normalsize
such that the covariance can be approximated 
and propagated 
over the horizon.
We obtain {\small$\mu^{d^j}, \Sigma^{x^j {d^j}}, \Sigma^{d^j}$} through a first-order Taylor approximation of the GP evaluated at its posterior mean \cite{kocijan_modelling_2016}, which provides a good trade-off between approximation accuracy and computational complexity \cite{hewing_cautious_2020}:

\vspace{-10pt}
\small
\begin{subequations} 
\begin{align}
    \mu_i^{d^j} & = \mu^{d^j}\left(\mu_i^{\mathbf{z}}\right), \vspace{-2pt} \\
    \Sigma_i^{x^j {d^j}}  & = \Sigma_i^{x^j}\Bigl( \nabla^\top_{x^j}\mu^{d^j}\left( \mu_i^{\mathbf{z}} \right) \Bigr)^\top, \vspace{-2pt} \\ 
    \Sigma_i^{d^j}  &= \Sigma^{d^j}\left( \mu_i^{\mathbf{z}} \right) + \nabla^\top_{x^j}\mu^{d^j}  \left( \mu_i^{\mathbf{z}} \right) \Sigma_i^{x^j}\Bigl(\nabla^\top_{x^j}\mu^{d^j} \left( \mu_i^{\mathbf{z}} \right)\Bigr)^\top ,
\end{align}
\end{subequations}%
\normalsize
where {\small$\mu^{\mathbf{z}}_i$} denotes the mean value of the test point {\small$\mathbf{z}_{i}$} at prediction step {\small$i$}, with {\small$\mu^{d^j}\left(\mathbf{z}\right)$} and {\small$\Sigma^{d^j}\left(\mathbf{z}\right)$} as in \eqref{eq:Sparse_GP}.
We propagate the mean and covariance of the posterior GP \eqref{eq:Sparse_GP} over the horizon through the following approximation \cite{hewing_cautious_2020}:

\vspace{-8pt}
\small
\begin{subequations}
\begin{align}
\mu_{i+1 \mid k}^{x^j} & = f^j \Bigl( \mu_{i \mid k}^{x^j} \Bigr) + B^j \mu^{d^j}\left( \mu^{\mathbf{z}}_{i \mid k}\right), \label{eq:Model_Mean} \\
\Sigma_{i+1 \mid k}^{x^j} & = \hspace{-0.5mm} \begin{bmatrix}
    \nabla^\top f^j\left( \mu_{i \mid k}^{x^j} \right) & \hspace{-2mm} B^j
\end{bmatrix} \Sigma^j_{i \mid k} \hspace{-0.5mm} \begin{bmatrix}
    \nabla^\top f^j\left( \mu_{i \mid k}^{x^j} \right) & \hspace{-2mm} B^j
\end{bmatrix}^{\top}
, \hspace{-0.5mm}
\label{eq:Model_Cov}
\end{align}
\label{eq:Sparse_GP_model}%
\end{subequations}
\normalsize
for {\small$i = 0, 1,\dots,N-1$},
where the initial prediction equals the current state {\small$\mu^{x^j}_{0 \mid k} = x^j_k$} with {\small${\Sigma^{x^j}_{ 0 \mid k} = \mathbf{0}}$}, and $\Sigma^j_{i \mid k}$ denotes the joint covariance matrix \eqref{eq:Joint_Dist}.
We approximate the full GP by a Sparse Pseudo-Input GP (SPGP) \cite{snelson_sparse_2005}, which exploits the Nyström projection to reduce the computational complexity of the GP by preconditioning the GP on \textit{inducing points}.
For details on SPGPs, we refer the reader to \cite{snelson_sparse_2005}.
We select the inducing points by taking equidistant samples from the predicted trajectory at time {\small$k-1$}, since the test points can be expected to be close to this trajectory \cite{hewing_cautious_2020}.

\subsection{Learning-based Dual Model Predictive Control}
In this section, we propose our dual MPC that uses the GP prediction model, presented in Sec. \ref{sec:GP-Model}, to jointly condition the interactions of the other agents with the control inputs of Agent $0$.
The dual control effect was first formalized in \cite{bar-shalom_dual_1974}, and according to \cite{mesbah_stochastic_2018} is defined as follows:

\begin{definition}
\label{def:DualControl}
    A control input is said to have the dual control effect if it can affect, with non-zero probability, at least one $r^{\text{th}}$-order central moment of a state variable {\small$(r \geq 2 )$}.
\end{definition}

We formulate a tractable dual MPC to devise a control policy for Agent $0$ and model the interactions $g^j$ by an SPGP:

\vspace{-9pt}
\small 
\begin{subequations}
\label{eq:Prim_MPC}
\begin{align}
\min _{U_k} \, J & \bigl(x^0_k, U_k \bigr) \\
\hspace{-2.8mm} \text { s.t. } x^0_{i+1 \mid k} &= f^0 \bigl(x^0_{i \mid k}, u^0_{i \mid k}\bigr), \quad i = 0, 1, \dots , N-1, \\
\mu^{x^j}_{i+1 \mid k} &= f^j \Bigl( \mu^{x^j}_{i \mid k} \Bigr) + B^j \mu^{d^j} \hspace{-1mm} \left( \mu^{\mathbf{z}}_{i \mid k} \right), \notag \\
i & = 0, 1, \dots, N-1, \quad j \in \mathcal{J}, \label{eq:MPC_Mean} \\
\Sigma_{i+1 \mid k}^{x^j} &= \begin{bmatrix} \nabla^\top \hspace{-0.5mm}f^j\left( \mu_{i \mid k}^{x^j} \right) & \hspace{-2mm} B^j \end{bmatrix} \Sigma^j_{i \mid k} \begin{bmatrix} \nabla^\top \hspace{-0.5mm} f^j\left( \mu_{i \mid k}^{x^j} \right) & \hspace{-2mm} B^j \end{bmatrix}^{\top}\hspace{-0.5mm}, \notag \\
i & = 0,1, \dots, N-1, \quad j \in \mathcal{J}, \label{eq:MPC_Cov} \\
u^0_{i \mid k} &\in \mathbb{U}^0, \quad i = 0,1, \dots , N-1, \\
x^0_{i \mid k} &\in \Tilde{{\mathbb{X}}}_{i \mid k}^j \Bigl(\mu_{i \mid k}^{x^j}, \Sigma^{x^j}_{i \mid k} \Bigr) , \quad i = 1,2, \dots , N, \quad j \in \mathcal{J}, \\
x^0_{i \mid k} & \in \mathbb{X}^0, \quad i = 1,2, \dots , N, \\
x^0_{0 \mid k} & = x^0_k, \\
\mu^{x^j}_{0 \mid k} &= x^j_k, \quad j \in \mathcal{J},\\
\Sigma^{x^j}_{0 \mid k} &= \boldsymbol{0},
\end{align}
\end{subequations}%
\normalsize
where the chance constraints from \eqref{eq:OCP_CC} are approximated by a general tightened constraint set {\small$\Tilde{\mathbb{X}}^j\bigl( \mu^{x^j}, \Sigma^{x^j}\bigr) \subseteq \mathbb{X}^j\bigl( \mu^{x^j}\bigr)$}, based on the covariance of the GP \cite{hewing_cautious_2020}. 
The specific form of constraint tightening may depend on the nature of the problem. An example hereof is provided in Sec. \ref{sec:Motion_Planning}.
The dual MPC exploits the joint conditioning of the predicted states of all the agents through the regression feature $\mathbf{z}$ of the GP. 
Consequently, the control sequence $U_k$ explicitly affects the constraint tightening $(\Tilde{\mathbb{X}}^j)$ through Bayesian inference.
The dual control effect provides an inherent degree of caution based on the correlation with past observations \cite{bar-shalom_dual_1974},
which is essential to a truly interaction-aware controller.
To this end, we make the following proposition. \vspace{7pt}

\begin{proposition}
The control input sequence of \text{Agent 0} {\small$(U_k^*)$} resulting from the locally optimal solution to the MPC problem \eqref{eq:Prim_MPC} affects both the trajectory of \text{Agent $0$} and the prediction of Agent {\small$j \in \mathcal{J}$} over the horizon. Moreover, {\small$U_k^*$} affects the predicted covariance of the state of Agent {\small$j \in \mathcal{J}$}, and, thus the control input {\small$u^{0*}_k$} has the dual control effect.

\end{proposition}
\begin{proof}
By construction, the control sequence resulting from \eqref{eq:Prim_MPC} {\small$U^*_{0:i-1 \mid k}$} determines the predicted future states of \text{Agent 0} {\small$x^0_{i \mid k}$}, and, consequently, the predicted test points {\small$\mathbf{z}_{i \mid k}$} of the GP. 
\vspace{-2pt}
Therefore, the control sequence {\small$U^*_{0:i-1 \mid k}$} affects the posterior mean {\small$\mu^{x^j}_{i \mid k}$} \eqref{eq:MPC_Mean} and posterior covariance {\small$\Sigma^{x^j}_{i \mid k}$} \eqref{eq:MPC_Cov} of the predicted states of Agent {\small$j \in \mathcal{J}$} at time {\small$k$} such that

\vspace{-5pt}
\small
\begin{equation}
    \Sigma^{x^j}_{i \mid k} \triangleq \mathbb{E}\bigl[\Sigma_{i}^{x^j} \mid \mathcal{D}_{k}, U^*_{0:i-1 \mid k} \bigr] \neq \mathbb{E}\bigl[ \Sigma_{i}^{x^j} \mid \mathcal{D}_{k} \bigr]
\end{equation}
\normalsize
for {\small$i = 1, 2, \dots, N$}. Hence, {\small$u^{0*}_k$} has the dual control effect.
\end{proof}

\section{Interaction-Aware Motion Planning}
\label{sec:Motion_Planning}
\subsection{Lane Merging Use Case}
In this section, we demonstrate the proposed GP-MPC method for the motion planning of an autonomous vehicle in a forced lane merging use case, and compare it against a stochastic baseline MPC.
This simple use case promotes the qualitative interpretability of GP-MPC.
Here, the target lane is occupied by a leading and a following vehicle,
as seen in \text{Fig. \ref{fig:Passive_Learning_Lapse}} in Sec. \ref{sec:Results}.
The Ego vehicle is controlled by the GP-MPC from \eqref{eq:Prim_MPC}.
The Follower's interactive driving policy follows a Merge-Reactive Intelligent Driver Model (MR-IDM) \cite{holley_mr-idm_2023}.
The Follower tries to close the gap with its Leader, hindering the Ego from merging.
We assume the target vehicles remain in their lane center.
The Ego, Follower, and Leader are denoted by {\small $j = 0,1,2$}, respectively.
For simplicity, we learn only the residual dynamics of the Follower {\small$(g^1)$}, assuming the Leader drives at constant speed {\small$(g^2 = 0)$}, and assume the disturbance processes {\small$w^j_k$} are zero.

\subsection{Vehicle Modeling}
All vehicles are modeled by a deterministic, kinematic bicycle model, which is sufficiently accurate for motion planning.
The state of a vehicle is described by the state vector
\linebreak
{\small$    x = \begin{bmatrix}
        X &
        Y &
        v &
        \psi &
        \delta
    \end{bmatrix}^\top \in \mathbb{R}^5,$}
where {\small$X$} and {\small$Y$} are the longitudinal and lateral positions of the rear axle, respectively. The longitudinal velocity of the rear axle is denoted by {\small$v$}, {\small$\psi$} is the heading angle of the vehicle, and {\small$\delta$} is the steering angle of the front axle.
The vehicle's acceleration $a$ and the steering rate $r$ are inputs to the system:
{\small${    u = \begin{bmatrix}
        a &
        r
    \end{bmatrix}^\top \in \mathbb{R}^2}$}.
The continuous-time dynamics of the kinematic bicycle model {\small${\dot{x}(t) = f_c\left(x(t), u \left(t \right)  \right)}$} 
are defined by the state-space equations:

\vspace{-10pt}
\small
\begin{subequations} 
\label{eq:CT_Dyn}
\begin{align}
\dot{X}(t) & = v(t) \cos \left(\psi\left(t\right)\right), \\
\dot{Y}(t) & = v(t) \sin \left(\psi\left(t\right)\right), \\
\dot{v}(t) & = a(t), \\
\dot{\psi}(t) & = \tfrac{v(t)}{l} \tan \left(\delta\left(t\right)\right), \\
\dot{\delta}(t) &= r(t),
\vspace{-15mm}
\end{align}
\end{subequations}
\normalsize
where $t \in \mathbb{R}$ denotes time and $l$ denotes the vehicle's wheelbase.
We discretize \eqref{eq:CT_Dyn} with the fourth-order Runge-Kutta method using a sampling time of ${T_s = 0.25 \text{ [s]}}$.

The Follower's MR-IDM policy \cite{holley_mr-idm_2023} is a driver model that considers agents merging into its lane.
This closed-loop policy {${a^1 = a_{\text{MR-IDM}} \bigl(\mathbf{x}, a^0, a^2\bigr)}$} couples the dynamics of the Ego and the Follower.
In order to utilize the framework outlined in Sec. \ref{sec:Method}, we reformulate this policy as a function of {$\mathbf{x}_k$} and {$u^0$} through \eqref{eq:Target_Pol}.
For simulation purposes, we assume the Follower has access to {$a^0$} and {$a^2=0$}.
The MR-IDM maps the relative longitudinal and lateral position to an effective distance gap for an IDM.
The vehicle inducing the largest deceleration determines the output of the MR-IDM.
For details on the MR-IDM we refer the reader to \cite{holley_mr-idm_2023}.

\subsection{Motion Planning Problem}
We incentivize the Ego to maintain its initial velocity {\small$v^0_0$} and stay in its lane, until the merge lane closes and the MPC decides that we should merge.
We use a smooth function {\small$m\left(X\right) : \mathbb{R} \rightarrow \mathbb{R}$} to describe the center of the merge lane as a function of the longitudinal position, which coincides with the center of the target lane after the merge lane has fully closed.
Accordingly, the Ego's reference $x^r$ is defined as:

\vspace{-3pt}
\small
\begin{equation}
\label{eq:state_ref}
    x^r\left( X \right) = \begin{bmatrix}
        0 & m(X) & v^0_0 & 0 & 0 
    \end{bmatrix}^\top.
    \vspace{-1mm}
\end{equation}
\normalsize
We employ the following objective function:

\vspace{-4pt}
\small
\begin{equation}
    \begin{aligned}
\label{eq:Obj_Fun}
J & \bigl(x^0_k, u^0_{k-1}, U_k \bigr) = \sum_{i=1}^{N} \norm{ x^0_{i \mid k} - x^r\left(X^0_{i \mid k}\right)}^2_Q \\
& + \sum_{i=0}^{N-1} \norm{u^0_{i\mid k}}^2_R + \norm{\Delta u^0_{i\mid k}}^2_{S},
\end{aligned}
\vspace{-5pt}
\end{equation}
\normalsize
where $\norm{z}^2_Q = z^\top Q z$ 
with positive (semi-) definite weighing matrices $Q, S \succeq 0$, $R \succ 0$ and $\Delta u_{i \mid k} = u_{i \mid k} - u_{i-1 \mid k}$ with $u_{-1 \mid k} = u_{k-1}$.
We limit the absolute acceleration and steering angle rate to $5$ $[\text{m/s}^2]$ and $5 \text{ [deg/s}] $ through 
$u \in \mathbb{U}^0$.
Constant state constraints partially govern road boundaries and limit the Ego's maximum velocity to $36$ [m/s], its absolute heading to $15 \text{ [deg]}$ and its absolute steering angle to $5 \text{ [deg]}$.
The right road boundary is dependent on its position and is governed through a separate constraint:

\vspace{-1mm}
\small
\begin{equation}
    h_r\bigl(x^0_{i \mid k} \bigr) = m\bigl(X^0_{i \mid k}\bigr) - Y^0_{i \mid k} + \tfrac{W - w}{2} \leq 0,
    \vspace{-1mm}
\end{equation}
\normalsize
where the road boundary is parallel to the center of the merge lane, $W$ and $w$ denote the vehicle and lane width.
Consequently, we have the following state constraint set:

\vspace{-10pt}
\small
\begin{equation}
\begin{aligned}
    \mathbb{X}^0 := \{x^0 \in \mathbb{R}^5 \mid x_{\text{min}} \leq x^0 \leq x_{\text{max}}, \, h_r\bigl(x^0 \bigr) \leq 0\}.
\end{aligned}
\vspace{-1mm}
\end{equation}
\normalsize

Similar to Zhu \textit{et al.} \cite{zhu_gaussian_2023}, we account for the uncertainty in the predicted position of the target vehicle by expanding the semi-axes of elliptical collision avoidance constraints.
Again, we assume no uncertainty in their lateral position.
By expanding the major semi-axis of the ellipse, we obtain a tractable, \textit{tightened} reformulation of the chance constraints:

\vspace{-8pt}
\small 
\begin{equation}
\label{eq:Expand_Ellipse}
\begin{aligned}
    & \Tilde{\mathbb{X}}_{i \mid k}^j \bigl(x_{i \mid k}^j \bigr) = \bigl\{ x^0_{i \mid k} \in \mathbb{R}^{n_x} \mid h_c\bigl(x^0_{i \mid k}, \mu^{x^j}_{i \mid k}, \Sigma^{X^j}_{i \mid k} \bigr) = \\
    & - \tfrac{\bigl( c^j_{x,i \mid k} - c^0_{x,i \mid k} \bigr) ^2}{ \left(\mathcal{E}_{c,A} + \sigma \sqrt{\Sigma^{X^j}_{i \mid k}} \right)^2 }
    - \tfrac{\bigl( c^j_{y,i \mid k} - c^0_{y,i \mid k} \bigr) ^2}{\mathcal{E}_{c,B}^2} + 1 \leq 0 \bigr\},
\end{aligned}
\vspace{-0.5mm}
\end{equation}
\normalsize
where {\small$c^j_{x,i \mid k}$} and {\small$ c^j_{y,i \mid k}$} denote the longitudinal and lateral component of the centroid, respectively, and {\small$\Sigma^{X^j}_{i \mid k}$} denotes the predicted covariance of the longitudinal position of \text{Agent $j$} at time step $k+i$.
The major and minor semi-axis of the ellipse are denoted by {\small$\mathcal{E}_{c,A}$} and {\small$\mathcal{E}_{c,B}$}, respectively.
The constraint tightening can be tuned with \text{\small$\sigma \in \mathbb{R}_{\geq 0}$}, where {\small$\sigma = 2$} is the amount of standard deviation that is accounted for.
In absence of formal feasibility guarantees of the MPC, the collision avoidance constraints are softened using an {\small$l_1$} penalty function
to retain feasibility.
Since we assume the Leader maintains its initial velocity, we set {\small $\Sigma^{x^2} = \boldsymbol{0}$}.

\label{sec:MPC}
\subsection{Constant Velocity Model Predictive Control}
The stochastic baseline MPC (CV-MPC) uses a constant velocity prediction model for both the Leader and Follower:

\vspace{-1mm}
\small
\begin{equation}    v^j_{i \mid k} = v^j_k, \quad \text{for } i = 0,1,\dots, N.
\vspace{-1mm}
\end{equation}
\normalsize
Hence, the nominal predictions $f^j$ of the target vehicles are:

\vspace{-10pt}
\small
\begin{equation}
\label{eq:Nominal_Pred}
    \hspace{-2.3mm}  \mu^{x^j}_{i+1 \mid k} = A \mu^{x^j}_{i\mid  k}, \quad \text{for } i = 0,1,\dots, N-1, \quad j = 1,2,
\vspace{-2mm}
\end{equation}
\normalsize
\vspace{-5pt}
\small \begin{equation}
\text{\normalsize where } A = \begin{bmatrix} 
        1 & 0 & T_s & 0 & 0 \\
        0 & 1 & 0   & 0 & 0 \\
        0 & 0 & 1   & 0 & 0 \\
        0 & 0 & 0   & 1 & 0 \\
        0 & 0 & 0   & 0 & 1 \\
    \end{bmatrix}.
\end{equation}
\normalsize
The CV-MPC accounts for the Follower's change in velocity {\small$\Delta v^1$} by a fixed covariance {\small${\sigma^2_{v^1} = 0.3}$} that we propagate as:

\small

\begin{equation}
    \Sigma_{i+1 \mid k}^{x^1} = A \Sigma_{i \mid k}^{x^1} A^\top + B^1 \sigma^2_{v^1} B^{1\top}, \quad 
\end{equation}
for $i = 0,1,\dots, N-1$, where
{\small$B^1 = \begin{bmatrix}
        0 & 0 & 1 & 0 & 0
        \end{bmatrix}^\top$}.
\normalsize

\subsection{Gaussian Process Model Predictive Control}
\normalsize

We use past observations ${\mathbf{y}}$ and predict the incremental velocity by the posterior of the GP through Bayesian inference.
Recall that the MPC only has access to the current state of other agents.
Therefore, the GP-MPC predicts the change of velocity by observing the difference between the Follower's actual velocity and the constant velocity prediction:

\vspace{-4pt}
\small
\begin{equation}
    g^1 \left(\mathbf{x}_k, \mathbf{u}_k \right) = y_k = v^1_{k+1} - v^1_{k} = \Delta v^1_k.
\end{equation}
\normalsize
We `select' specific states of interest by applying a linear map $C$ to $\mathbf{z}$, such that we infer these predictions using the velocities and relative positions of the different vehicles:

\vspace{-5pt}
\small 
\begin{equation}
\label{eq:test_point_AV}
    C \hspace{0.5mm} \mathbf{z} = \hspace{-0.5mm} \begin{bmatrix}
        v^0,
        v^1,
        v^2,
        \bigl( X^1 - X^0 \bigr),
        \bigl( X^1 - X^2 \bigr),
        \bigl( Y^1 - Y^0 \bigr)
    \end{bmatrix}^\top \hspace{-2mm}. \hspace{-1mm}
\end{equation}
\normalsize
The relative lateral position of the Follower and Leader is assumed to be zero and is omitted from training point $\mathbf{z}$.

We employ the frequently used squared exponential kernel and map the regression features through $C$, such that:

\vspace{-5pt}
\small
\begin{equation}
\label{eq:SE_ker}
    \hspace{-3mm} k\left(\mathbf{z}, \mathbf{z}'\right) = \sigma^2_d \exp \bigl( - \tfrac{1}{2} \left(C\mathbf{z} \hspace{-0.3mm} -\hspace{-0.3mm} C\mathbf{z}' \right)^\top \hspace{-0.4mm} L_d^{-2} \left(C\mathbf{z} \hspace{-0.3mm}-\hspace{-0.3mm} C\mathbf{z}'\right) \bigr) , \hspace{-2mm}
\end{equation}
\normalsize
where $\sigma^2_d$ is the prior covariance and {\small $L_d$} is the length-scale matrix.
If sufficient training data is available, these hyper-parameters can be inferred using log-likelihood optimization \cite{rasmussen_gaussian_2006}.
Here, we use fixed hyper-parameters for the kernel with {\small$ L_d = \text{diag}\left(3,3,3,17,17,5\right)$}.
The prior covariance is equal to the covariance of the baseline prediction model {\small$\sigma^2_d = 0.3$}, making the CV-MPC equivalent to the GP-MPC prior to inference. Since we do not consider any noise, we set {\small$\sigma^2_w = 0$}.

\section{Results}
\label{sec:Results}
We compare our proposed GP learning-based dual MPC (GP-MPC) with a stochastic baseline MPC that uses a constant velocity prediction (CV-MPC)
for the Follower, and evaluate their performance on the lane merging use case detailed in Sec. \ref{sec:Motion_Planning}.
We test the GP-MPC with and without pre-training. 
Here, we refer to pre-training as the availability of past measurement data at the start of the experiment.

\subsection{Numerical Case Study}
We perform 51 Monte Carlo simulations to assess the prediction capabilities of the GP-MPC and CV-MPC.
As initial conditions, we take {\small$ X^1_0 = -75 \text{ [m]}, v^1_0 = 31 \text{ [m/s]} $} for the Follower, and {\small$ X^2_0 = 0 \text{ [m]}, v^2_0 = 25 \text{ [m/s]}$} for the Leader to ensure a sufficiently challenging scenario. The initial position of the Ego is equidistantly sampled from {\small$ X^0_0 \in \left[ -100, -75 \right]$} with {\small$ v^0_0 = 31 \text{ [m/s]}$}.
Computations are performed on an Intel i7 with 32 GB RAM, in MATLAB
using CasADi \cite{andersson_casadi_2019}, with the IPOPT \cite{wachter_implementation_2006} and MA57 \cite{hsl_coin-hsl_2023} solvers. We use a custom GP implementation in CasADi to compute the joint covariance. We take a horizon of {\small$N=12$}, {\small$M=4$} inducing points for the SPGP and a sampling time of {\small$T_s = 0.25$} [s].
At each time {\small$k+1$}, we append the training data: {\small$\mathcal{D}_{k+1} = \mathcal{D}_k \cup (\mathbf{z}_{k}, \mathbf{y}_k)$}. Although we currently do not discard any data points, selection algorithms can help reduce the size of the training set \cite{kabzan_learning-based_2019}.
To assess the generalizability to unseen behavior without pre-training, we take {\small$ \mathcal{D}_0 = \emptyset$}.
With pre-training, we initialize the GP with {\small$n_D = 80$} observations from a single closed-loop experiment with {\small$ X^0_0 = -85 \text{ [m]}$}.

\begin{table}[t]
\centering
\caption{Results of the Monte Carlo simulations.}
\vspace{-2mm}
\begin{tabular}{l c c c c} \hline \hline \\[-2.5mm]
Algorithm  & Success & $\overline{\|e\|}$  $[\text{m/s}]$ & $\overline{T}_c$ $[\text{ms}]$ & $T_c < T_s$ $[\%]$
\\[0.5mm] \hline
CV-MPC     & 20/51   & $0.679$  & $50$  & $100$ \\ \hline
GP-MPC     & 33/51   & $0.645$  & $128$ & $95.6$ \\ \hline
GP-MPC with \hspace{-2mm} \\ pre-training   & 35/51   & $0.440$    & $206$ & $69.0$ \\ \hline
\end{tabular}
\label{tab:Results_MC}
\vspace{-10pt}
\end{table}

Table \ref{tab:Results_MC} summarizes the results of the Monte Carlo trials.
In this study, we consider a merge successful if the Ego merges between the Follower and Leader.
The prediction error {\small ${\|e\left(k\right)\| := \frac{1}{N} \| \{v_{i \mid k} - v_{k+i} \}_{i=1}^N \|_{1} }$} is the normalized norm of the difference between the predicted trajectory at time $k$ and the realized trajectory from time $k$ to $k+N$.
The average prediction error over all time steps and all trials {\small $\overline{\|e\|}$} quantifies the prediction quality of the CV-MPC and GP-MPC.
The average solve time {\small$\Bar{T}_c$} and the percentage of iterations within sampling time {\small$T_c < T_s$} show that the algorithm can be run (near-)realtime with these amounts of training data.

\subsection{Discussion}
\begin{figure}[t]
    \centering
    \includegraphics[width = \columnwidth]{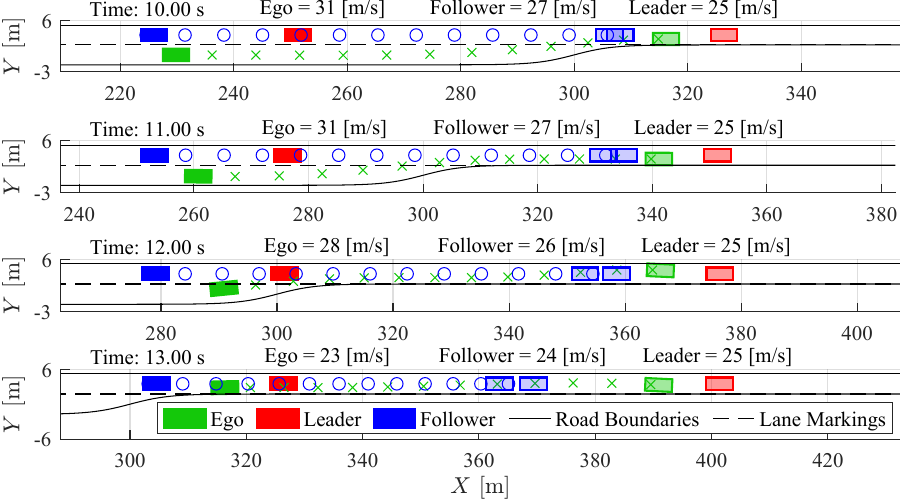}
    \vspace{-20pt}
    \caption{Time-lapse of a lane merge by the GP-MPC. The transparent rectangles show the final predicted positions of vehicles, while for the Follower two rectangles indicate \text{$2\sigma$-bounds} on its final predicted position.}
    \label{fig:Passive_Learning_Lapse}
    \vspace{-12pt}
\end{figure}

The GP-MPC jointly optimizes the predictions and the associated uncertainty through the dual control effect, as seen in
\text{Fig. \ref{fig:Passive_Learning_Lapse}}, which shows a successful merge by the GP-MPC without pre-training.
The approximate covariance of the GP is exploited to account for uncertain interactions providing inherent caution through dual MPC.
As the Ego is agnostic to the Follower's response to the merging action, it takes extra caution when merging by further tightening the constraints.

\begin{figure}[t]
  \centering
  \includegraphics[width = \columnwidth]{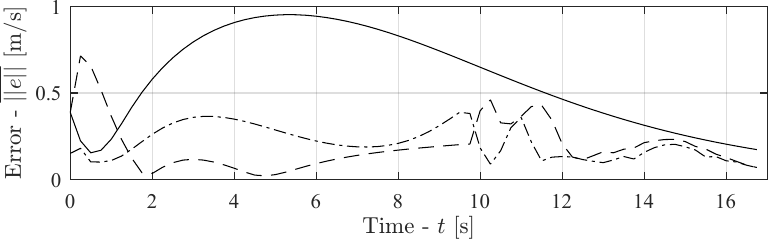}  \label{fig:a} \\ \vspace{-8pt}
  \includegraphics[width = \columnwidth]{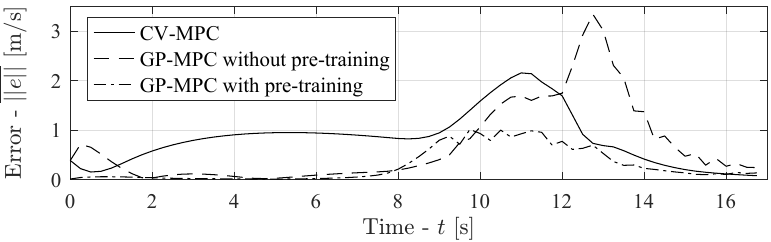} \label{fig:b}
  \vspace{-15pt}
  \caption{Average prediction error of Follower's velocity over the prediction horizon in scenarios in which all planners merge behind (average over 16 experiments, shown in the top figure) and merge in between (average over 20 experiments, shown in the bottom figure).} \label{fig:Pred_error_between}
  \vspace{-10pt}
\end{figure}

Both the CV-MPC and the GP-MPC without pre-training do not capture the high frequent dynamics during the merge, as seen in Fig. \ref{fig:Pred_error_between}.
During the merge, approximately after {\small$t=10$} [s], the GP-MPC without pre-training has an increased prediction error as it has not observed these high frequent dynamics.
Still, the GP-MPC has superior prediction quality prior to merging and is able to identify a safe gap to merge in between.
The GP-MPC with pre-training exploits past observations from a similar scenario, leveraging Bayesian inference and the dual control effect to exercise less caution and shows increased prediction quality.
When compared to CV-MPC, the online learning GP-MPC sees a $25 \%$ to $29 \%$ increase in successful merges over all experiments without and with pre-training, respectively, cf. \text{Table \ref{tab:Results_MC}}.
While GPs are attractive for Bayesian inference, the computational complexity increases with the number of training data. This stresses the importance of online learning in dynamic environments, such that the MPC does not have to rely on large sets of training data, but rather adapt online.
Simulation results show that our GP-MPC can successfully predict the Follower's motion from various initial conditions by only relying on online training data, demonstrating its generalizability beyond a priori available training data.

\section{Conclusions}
We presented a Gaussian process-based online learning MPC that uses the dual control effect to learn and predict the dynamics of other agents while accounting for their joint uncertainty.
We validated our GP-MPC in a numerical case study on automated lane merging, by leveraging dual MPC to jointly optimize the motion of the ego vehicle and the predictions of a target vehicle through an interactive GP model.
By explicitly considering the effect of control inputs on the covariance of predictions of other agents, our GP-MPC can leverage this effect in the tightening of constraints.
Through online learning, the GP-MPC can adapt to uncertain and unseen behavior of other agents, showing strong potential for scenarios that extend beyond a priori available training data sets.

The dual control effect is an important notion that can be leveraged by the MPC.
While this paper focuses on a proof of concept of the proposed method, the next research steps for learning-based GP-MPC include:
(i) incorporating the probing effect of dual control, i.e., active learning,
\text{(ii) predicting} the residual dynamics of multiple agents, and (iii) providing performance and recursive feasibility guarantees of the MPC. In the context of autonomous driving, planned future work includes validation of the prediction model, the learning of various driving behaviors, and experimental validation.






\bibliographystyle{IEEEtran}
\bibliography{Ref}

\end{document}